\documentclass[12pt]{article}

\usepackage{latexsym}
\usepackage{amsfonts}
\usepackage{amssymb}
\usepackage{amsmath}
\usepackage{mathrsfs}
\usepackage{hyperref}

\newcommand{\be}{\begin{equation}}
\newcommand{\ee}{\end{equation}}
\newcommand{\bea}{\begin{eqnarray}}
\newcommand{\eea}{\end{eqnarray}}
\newcommand{\bean}{\begin{eqnarray*}}
\newcommand{\eean}{\end{eqnarray*}}
\newcommand{\brray}{\begin{array}}
\newcommand{\erray}{\end{array}}

\newtheorem{dfn}{Definition}[section]
\newtheorem{thm}[dfn]{Theorem}
\newtheorem{lmma}[dfn]{Lemma}
\newtheorem{ppsn}[dfn]{Proposition}
\newtheorem{crlre}[dfn]{Corollary}
\newtheorem{xmpl}[dfn]{Example}
\newtheorem{rmrk}[dfn]{Remark}

\newcommand{\bdfn}{\begin{dfn}\rm}
\newcommand{\bthm}{\begin{thm}}
\newcommand{\blmma}{\begin{lmma}}
\newcommand{\bppsn}{\begin{ppsn}}
\newcommand{\bcrlre}{\begin{crlre}}
\newcommand{\bxmpl}{\begin{xmpl}}
\newcommand{\brmrk}{\begin{rmrk}\rm}

\newcommand{\edfn}{\end{dfn}}
\newcommand{\ethm}{\end{thm}}
\newcommand{\elmma}{\end{lmma}}
\newcommand{\eppsn}{\end{ppsn}}
\newcommand{\ecrlre}{\end{crlre}}
\newcommand{\exmpl}{\end{xmpl}}
\newcommand{\ermrk}{\end{rmrk}}

\newcommand{\clh}{\mathcal{H}}

\newcommand{\clk}{\mathcal{K}}

\author{S.P. Murugan and S. Sundar}
\title{On the existence of $E_{0}$-semigroups - the multiparameter case}
\date{}

\begin{document}
\maketitle
\begin{abstract}
Let $P \subset \mathbb{R}^{d}$ be a closed convex cone. Assume that $P$ is pointed, i.e. the intersection $P \cap -P=\{0\}$ and $P$ is spanning, i.e. $P-P=\mathbb{R}^{d}$. Denote the interior of $P$ by $\Omega$. 
Let $E$ be a product system  over  $\Omega$. We show that there exists an infinite dimensional  separable Hilbert space $\clh$
and a semigroup $\alpha:=\{\alpha_x\}_{x \in P}$ of unital normal $*$-endomorphisms of $B(\clh)$ such that $E$
is isomorphic to the product system associated to $\alpha$. 
\end{abstract}
\noindent {\bf AMS Classification No. :} {Primary 46L55; Secondary 46L99.}  \\
{\textbf{Keywords :}}$ E_{0}$-semigroups, Essential representations, Product systems.

\section{Introduction}
Let $P \subset \mathbb{R}^{d}$ be a closed convex cone which is spanning, i.e. $P-P=\mathbb{R}^{d}$ and pointed, i.e. $P \cap -P=\{0\}$. Denote the interior of $P$ by $\Omega$. 
Then $\Omega$ is an ideal in $P$ in the sense that $\Omega + P \subset \Omega$. Moreover $\Omega \cap -\Omega = \emptyset$. Also $\Omega$ is dense in $P$. We reserve the above notation for the rest of 
this paper. All the Hilbert spaces that we consider are separable and are over the field of complex numbers. For a Hilbert space $\clh$, $B(\clh)$ denotes the algebra of bounded operators on $\clh$. Let $\clh_1$ and $\clh_2$ be Hilbert spaces and $U:\clh_1 \to \clh_2$ be a unitary. We denote the map $B(\clh_1) \ni T \to UTU^{*} \in B(\clh_2)$ by $Ad(U)$.

Let $\clh$ be an infinite dimensional separable Hilbert space. 
By an $E_{0}$-semigroup over $P$ on $B(\clh)$, we mean a family $\alpha:=\{\alpha_{x}\}_{x \in P}$ of normal $*$-endomorphisms of $B(\clh)$ such that 
\begin{enumerate}
\item[(1)] for $x,y \in P$, $\alpha_{x} \circ \alpha_{y}=\alpha_{x+y}$, 
\item[(2)] for $x \in P$, $\alpha_{x}$ is unital, i.e. $\alpha_{x}(1)=1$, and
\item[(3)] for $A \in B(\clh)$ and $\xi,\eta \in \clh$, the map $P \ni x \to \langle \alpha_{x}(A)\xi|\eta \rangle \in \mathbb{C}$ is continuous.
\end{enumerate}
When $P=[0,\infty)$, we recover the classical theory of $1$-parameter $E_0$-semigroups which has a rich history. We refer the reader to the monograph \cite{Arveson} and the references therein for its 
literature.  In comparison to the literature on $1$-parameter $E_0$-semigroups, the literature on $E_0$-semigroups parametrised by semigroups other than $\mathbb{N}$ and $\mathbb{R}_{+}$ are few. The notable ones are \cite{Solel}, \cite{Shalit_2008},  \cite{Shalit}, \cite{Shalit_Solel}, \cite{Shalit_2011}, \cite{Hirshberg_Daniel} and \cite{Vernik}. We should also mention here that there is a good amount of progress achieved in developing the theory of  one parameter $E_0$-semigroups on operator algebras other than the algebra of bounded operators on a Hilbert space. The important papers in this direction are \cite{Skeide_Bhat},  \cite{Solel_Muhly}, \cite{Skeide2},  \cite{Srinivasan_Oliver1}, \cite{Skeide4} and  \cite{Srinivasan_Oliver2}. It would be interesting to develop the multi-parameter $E_0$-semigroup theory in this context. In this paper, we consider only $E_0$-semigroups over $P$ on $B(\clh)$. To avoid repetition of the phrase, ``$E_0$-semigroups over $P$", we will call an $E_0$-semigroup over $P$ simply an $E_0$-semigroup.

Let $\alpha:=\{\alpha_{x}\}_{x \in P}$ be an $E_{0}$-semigroup on $B(\clh)$. A strongly continuous family of unitaries $U:=\{U_{x}\}_{x \in P}$ in $B(\clh)$ is called an $\alpha$-cocycle if $U_{x}\alpha_{x}(U_y)=U_{x+y}$ for $x,y \in P$.  Let $U:=\{U_{x}\}_{x \in P}$ be an $\alpha$-cocycle. Then $\{Ad(U_x)\circ \alpha_{x}\}_{x \in P}$ is an $E_{0}$-semigroup. Such an $E_{0}$-semigroup is called a cocycle perturbation of $\alpha$. Let $\beta:=\{\beta_{x}\}_{x \in P}$ be an $E_{0}$-semigroup on $B(\clk)$ where $\clk$ is an infinite dimensional separable Hilbert space. We say $\alpha$ and $\beta$ are cocycle conjugate if there exists a unitary $U:\clh \to \clk$ such that $\{Ad(U^{*})\circ \beta_{x} \circ Ad(U)\}_{x \in P}$ is a cocycle perturbation of $\alpha$. It is routine to see that cocycle conjugacy is an equivalence relation. 

The main problem in the theory of $E_0$-semigroups is to classify them up to cocycle conjugacy. Arveson found a complete invariant called the product system associated to an $E_{0}$-semigroup when $P=[0,\infty)$.  Let us explain the notion of a product system associated to an $E_{0}$-semigroup. Let $\alpha:=\{\alpha_{x}\}_{x \in P}$ be an $E_{0}$-semigroup on $B(\clh)$. 
For $x \in \Omega$, let \[
E(x):=\{T \in B(\clh): \alpha_{x}(A)T=TA ~\textrm{for~}A \in B(\clh)\}.\]
We endow $B(\clh)$ with the $\sigma$-algebra generated by  weakly closed subsets of $B(\clh)$. The product $\Omega \times B(\clh)$ is given the product $\sigma$-algebra where  the $\sigma$-algebra on $\Omega$ is the Borel $\sigma$-algebra.   

Let \[E:=\{(x,T): x \in \Omega, T \in E(x)\}\] and let $p:E \to \Omega$ be the first projection.  For $x \in \Omega$, we identify $p^{-1}(x)$ with $E(x)$. Then we have the following.
\begin{enumerate}
\item[(1)] The set $E$ is a measurable subset of the standard Borel space $\Omega \times B(\clh)$.
\item[(2)] Fix $x \in \Omega$. Then for $S,T \in E(x)$, $T^{*}S$ is a scalar which we denote by $\langle S|T \rangle$. With respect to the inner product $\langle~|~\rangle$, the vector space $E(x)$ is a separable Hilbert space.
\item[(3)] For $x,y \in \Omega$, the linear span of $\{ST: S \in E(x), T \in E(y)\}$ is dense in $E(x+y)$.
\item[(4)] There exists a separable Hilbert space $\clh_0$ such that the following holds.  For $x \in \Omega$, there exists a unitary operator $\theta_{x}:E(x) \to \clh_0$ such that the map \[E \ni (x,T) \to (x,\theta_{x}(T))\in \Omega \times \clh_0\] is a Borel isomorphism where the measurable structure on $\Omega \times \clh_0$ is the one induced by the product topology on $\Omega \times \clh_0$. Here $\clh_0$ is given the norm topology.
\item[(5)] The set $E$ has an associative multiplication given by the formula \[(x,S).(y,T)=(x+y,ST)\] for $x,y \in \Omega$ and $(S,T) \in E(x) \times E(y)$. 
\end{enumerate}
The set $E$ together with the above structures is called the product system associated to $\alpha$. 
For a proof of the above facts, we refer the reader to \cite{Murugan_Sundar}. It is easy to see that $(1)$, $(2)$ and $(3)$ holds for any $E_0$-semigroup. Strictly speaking, $(4)$ is established in \cite{Murugan_Sundar} under the assumption that $\alpha_{a}$ is not onto for every $a \in \Omega$.  If $\alpha_{a}$ is onto for some $a \in \Omega$ then by Lemma 2.5 of \cite{Murugan_Sundar}, it follows that $\alpha_{x}$ is an automorphism for every $x \in P$. Then $(4)$ is an easy consequence of Arveson-Wigner's  Theorem (Theorem 3.6 of \cite{Anbu}). The statement of Arveson-Wigner's Theorem is as follows : Let $\alpha:=\{\alpha_{x}\}_{x \in P}$ be an $E_0$-semigroup on $B(\clh)$ such that for every $x \in P$, $\alpha_{x}$ is an automorphism. Then there exists a strongly continuous family of unitaries $U:=\{U_{x}\}_{x \in P}$ and a strictly upper triangular $d \times d$ real matrix $A$ such that $\alpha_{x}(T)=U_{x}TU_{x}^{*}$ for $T \in B(\clh)$ and $U_{x}U_{y}=e^{i\langle Ax|y \rangle}U_{x+y}$ for $x, y \in P$. Though not needed in this paper, we must mention here that the $E_0$-semigroups on $B(\clh)$ consisting of automorphisms, up to cocycle conjugacy, are in one-one correspondence with the set of strictly upper triangular $d \times d$ real  matrices.  A proof of this can be found in \cite{Anbu}. 

The authors, imitating the $1$-dimensional proof of Arveson,  proved in \cite{Murugan_Sundar} that two $E_0$-semigroups are cocycle conjugate if and only if the associated product systems are isomorphic. Just like in the $1$-dimensional case, we can define a product system over $\Omega$ abstractly, i.e. an abstract product system over $\Omega$ is a standard Borel space with structures reflecting the structures of a product system associated to an $E_0$-semigroup. Then it is natural to ask whether an abstract product system over $\Omega$ is isomorphic to a product system associated to an $E_0$-semigroup. The goal of this paper is to answer this question in the affirmative. 

The case when $P=[0,\infty)$ was first settled by Arveson using the machinery of the spectral $C^{*}$-algebra of a product system. The proof is technical and long. Later Skeide in \cite{Skeide} found a simpler proof.  Arveson himself found a simpler proof in \cite{Arv}. We should also mention that Skeide proved in \cite{Skeide1} that the $E_0$-semigroups arising out of  his construction and Arveson's construction are in fact conjugate and not just cocycle conjugate. Here we imitate Arveson's proof in \cite{Arv}. The authors in \cite{Murugan_Sundar_discrete} have done a similar analysis in the discrete setting for a finitely generated subsemigroup of $\mathbb{Z}^{d}$.  
We must mention here that this paper (and also its title) is heavily inspired by \cite{Arv}. 

We would like to thank  R. Srinivasan of Chennai Mathematical Institute, the Ph.d advisor of the first author, for introducing us to the beautiful world of $E_0$-semigroups and product systems.
We thank  Partha Sarathi Chakraborty, Institute of Mathematical Sciences for his geometric insight which turned out to be crucial in proving Lemma \ref{shift}. We thank the anonymous referee for his/her careful reading of the earlier version of the paper and his/her suggestions which  has significantly improved the readability of this paper. 

\section{Preliminaries}

Let $P_1 \subset \mathbb{R}^{d_1}$ and $P_2 \subset \mathbb{R}^{d_2}$ be closed convex cones. A map $\phi:P_1 \to P_2$ is called an isomorphism if $\phi$ is a bijection,  $\phi(x+y)=\phi(x)+\phi(y)$, $\phi(\lambda x)=\lambda \phi(x)$ for $x,y \in P_1$ and $\lambda \geq 0$. The cones $P_1$ and $P_2$ are said to be isomorphic if there exists a map $\phi:P_1 \to P_2$ such that $\phi$ is an isomorphism. Let $<P_i>=P_i-P_i$ for $i=1,2$. It is routine to check the following two facts. 
\begin{itemize}
\item For $i \in \{1,2\}$, $<P_i>$ is a vector subspace of $\mathbb{R}^{d_i}$.
\item Suppose that $\phi:P_1 \to P_2$ is an isomorphism. Then there exists a unique linear isomorphism from $<P_1>$ to $<P_2>$ which extends $\phi$. 

\end{itemize}

Although we do not deal with any particular example of a cone, let us give a few examples of cones that the reader should keep in mind. 

\begin{enumerate}
\item[(1)]Let $\mathbb{R}_{+}^{d}:=\{(x_1,x_2,\cdots,x_d):~\textrm{for every $i$,~}x_i \geq 0\}$. Clearly $\mathbb{R}_{+}^{d}$ is a closed convex cone which is both spanning and pointed.
\item[(2)] Let $V$ be a finite dimensional real inner product space. Denote the space of symmetric linear operators on $V$ by $\mathcal{S}(V)$ and the cone of positive operators on $V$ by
               $\mathcal{P}(V)$. Then it is clear that $\mathcal{P}(V)$ is a closed convex cone in $\mathcal{S}(V)$ which is both spanning and pointed.
 \item[(3)] Let $d \geq 3$ and $Q:=\{(x_1,x_2,\cdots,x_d):x_d \geq 0, x_d^{2} \geq \sum_{k=1}^{d-1}x_{k}^{2}\}$. Then $Q$ is a closed convex cone in $\mathbb{R}^{d}$ which is both pointed and spanning. The cone $Q$ is often called the forward light cone in $\mathbb{R}^{d}$.
\end{enumerate}

 The cones described in $(1)$, $(2)$ and $(3)$ are examples of what are called symmetric cones. There is a well known bijective correspondence between symmetric cones and Euclidean Jordan algebras from which it follows that the cones described in $(1)$, $(2)$ and $(3)$ are mutually non-isomorphic. We will not enter into a long digression trying to distinguish the cones in $(1)$, $(2)$ and $(3)$ by elementary means as it is not relevant to the rest of this article. To know more about the exact relationship between Jordan algebras and symmetric cones, we refer the reader  to the first five chapters of the monograph \cite{Faraut}. Of course, there are other possibilites like polyhedral cones whose classification we are not aware  of if one such exists. We satisfy ourself with  a reasoning  that demonstrates that the cone of positive operators on a finite dimensional inner product space is not isomorphic to $\mathbb{R}_{+}^{k}$ for any $k$ when the dimension of the underlying vector space is atleast 2.

Let $P \subset \mathbb{R}^{d}$ be a closed convex cone and $F \subset P$ be a closed subsemigroup of $P$  containing the origin $0$. We say that $F$ is a face of $P$ if $F$ is closed and if $y+z \in F$ with $y,z \in P$ then $y,z \in F$. It is routine to see that a face of $P$ is always a closed convex cone. Let $F \subset P$ be a face. Then $F$ is said to be one dimensional if the linear span of $F$ i.e. $<F>=F-F$ is one dimensional. Note the following.

\begin{enumerate}
\item[(a)] Let $V$ be a finite dimensional real inner product space of dimension atleast 2. Denote the linear space of symmetric linear operators on $V$ by $\mathcal{S}(V)$ and the cone of positive operators on $V$ by $\mathcal{P}(V)$. Denote the set of one dimensional faces of $\mathcal{P}(V)$ by $\mathcal{F}$ and the set of rank one projections on $V$ by $\mathcal{Q}$. It follows immediately from the spectral theorem that the map $\mathcal{Q} \ni Q \to \{\lambda Q: \lambda \geq 0\} \in \mathcal{F}$ is a bijection.
This implies in particular that the set of one dimensional faces is uncountable.
\item[(b)] On the other hand, the set of one dimensional faces of $\mathbb{R}_{+}^{k}$ has cardinality $k$.
\end{enumerate}
Hence, it follows that the cone of positive operators on a finite dimensional inner product space is not isomorphic to $\mathbb{R}_{+}^{k}$ for any $k \geq 1$ when the dimension of the underlying vector space is atleast 2. 

Henceforth for the rest of this paper, we reserve the letter $P$ to denote a closed convex cone in $\mathbb{R}^{d}$ which is spanning and pointed and the letter $\Omega$ to denote the interior of $P$. 
 We now introduce the  notion of an abstract product system over $\Omega$. We imitate Arveson's definition in \cite{Arveson} (Page 68, Definition 3.1.1.).

\begin{dfn}
\label{abstract}
By an abstract product system over $\Omega$, we mean a standard Borel space $E$ together with a measurable surjection $p:E \to \Omega$ such that the following holds.
\begin{enumerate}
\item[(1)] For $x \in \Omega$, $E(x):=p^{-1}(x)$ is a non-zero separable Hilbert space.
\item[(2)] There exists an associative multiplication $E \times E \ni (u,v) \to uv \in E$ such that $p(uv)=p(u)+p(v)$ for $u,v \in E$. Also, the multiplication $E \times E  \ni (u,v) \to uv \in E$ is measurable. 
\item[(3)] Let $x,y \in \Omega$ be given. Then there exists a unitary $u_{x,y}:E(x)\otimes E(y) \to E(x+y)$ such that $u_{x,y}(u \otimes v)=uv$ for $(u,v) \in E(x) \times E(y)$. 
\item[(4)] Let $\Delta:=\{(u,v) \in E\times E: p(u)=p(v)\}$. The maps $\Delta \ni (u,v) \to u+v \in E$ and $\Delta\ni (u,v) \to \langle u|v \rangle \in \mathbb{C}$ are measurable. 
\item[(5)] The map $\mathbb{C} \times E \ni (\lambda,u) \to \lambda u \in E$ is measurable.
\item[(6)] As a measurable field of Hilbert spaces, $E$ is trivial by which we mean the following. There exists a separable Hilbert space $\clh_0$ and a Borel isomorphism $\theta:E \to \Omega \times \clh_0$ such that $\pi_1 \circ \theta=p$ and for every $x \in \Omega$, the map $\pi_2 \circ \theta:E(x) \to \clh_0$ is a unitary. Here by $\pi_1$ and $\pi_2$, we mean the first and second projections from $\Omega \times \clh_0$ onto $\Omega$ and $\clh_0$ respectively. We shall abbreviate the above by saying that as a measurable field of Hilbert spaces, $E$ is isomorphic to $\Omega \times \clh_0$ and write $E \simeq \Omega \times \clh_0$.
\end{enumerate}
\end{dfn}
We usually suppress the surjection $p$ and simply write $E$ as $\displaystyle E=\coprod_{x \in \Omega}E(x)$.

\begin{rmrk}
Let $\clh$ be an infinite dimensional separable Hilbert space and $\alpha:=\{\alpha_{x}\}_{x \in P}$ be an $E_0$-semigroup. Then the product system associated to $\alpha$ is an abstract product system.

\end{rmrk}

\begin{dfn}
 Let $\displaystyle E:=\coprod_{x \in \Omega}E(x)$ and $\displaystyle  F:=\coprod_{x \in \Omega}F(x)$ be abstract product systems over $\Omega$. We say that $E$ is isomorphic to $F$ if for every $x \in \Omega$, there exists a unitary map $\theta_{x}:E(x) \to F(x)$ such that 
\begin{enumerate}
\item[(1)] for $x,y \in \Omega$ and $(u,v) \in E(x) \times E(y)$, $\theta_{x+y}(uv)=\theta_{x}(u)\theta_{y}(v)$, and
\item[(2)] the map $\displaystyle \theta:=\coprod_{x \in \Omega}\theta_{x}:E \to F$ is a Borel isomorphism. 
\end{enumerate}
Since $E$ and $F$ are standard Borel spaces, it is enough to require that the map $\theta$ in $(2)$ is $1$-$1$, onto and measurable (See Page 70, Theorem 3.3.2. of \cite{Arveson_invitation}). 
\end{dfn}

Let us make a few preliminary observations regarding the dimension of the fibres of an abstract product system. We will also drop the adjective ``abstract" and simply call an abstract product system over $\Omega$ a product system over $\Omega$.

\begin{lmma}
Let $\displaystyle E:=\coprod_{x \in \Omega}E(x)$ be a product system over $\Omega$. For $x \in \Omega$, let $d(x)$ be the dimension of $E(x)$. Then either 
$d(x)=1$ for all $x \in \Omega$ or $d(x)=\infty$ for all $x \in \Omega$. 
\end{lmma}
\textit{Proof.} Note that $d(x) \neq 0$ for every $x \in \Omega$. Condition $(6)$ of Definition \ref{abstract} implies that $d(x)=d(y)$ for every $x,y \in \Omega$. 
Thus it is enough to show that $d(a)=1$ or $d(a)=\infty$ for some $a \in \Omega$. Let $a \in \Omega$ be given. Then $E(a) \otimes E(a) \cong E(2a)$. Hence $d(a)=d(2a)=d(a)^{2}$. This implies $d(a)=1$ or $d(a)=\infty$. This completes the proof. \hfill $\Box$

The goal of this paper is to show that every abstract product system over $\Omega$ is isomorphic to the product system associated to an $E_0$-semigroup on $B(\clh)$ where $\clh$ is an infinite dimensional separable Hilbert space. First,  we dispose of the case when the fibres are $1$-dimensional. 

\begin{ppsn}
Let $\displaystyle E:=\coprod_{x \in \Omega}E(x)$ be a product system over $\Omega$ such that  the fibre $E(x)$ is one-dimensional for every $x \in \Omega$. Then there exists an $E_0$-semigroup  $\alpha:=\{\alpha_{x}\}_{x \in P}$ on $B(L^{2}(\mathbb{R}^{d}))$ such that for  $x \in P$, $\alpha_{x}$ is an automorphism and the product system associated to $\alpha$ is isomorphic to $E$.
\end{ppsn}
\textit{Proof.} Since $E \simeq \Omega \times \mathbb{C}$, it follows that there exists a measurable section $e:\Omega \to E$ such that for $x \in \Omega$, $||e(x)||=1$ and $E(x)$ is spanned by $e(x)$. Let $x,y \in \Omega$ be given. Then there exists a unique scalar denoted $\omega(x,y) \in \mathbb{T}$ such that $e(x)e(y)=\omega(x,y)e(x+y)$. Observe that for $x,y \in \Omega$, $\omega(x,y)=\langle e(x)e(y)|e(x+y) \rangle$. This implies that the function $\Omega \times \Omega \ni (x,y) \to \omega(x,y) \in \mathbb{T}$ is measurable. The associativity of the multiplication of the product system implies that $\omega$ is a multiplier on $\Omega$, i.e. for $x,y,z \in \Omega$, 
\[
\omega(x,y)\omega(x+y,z)=\omega(x,y+z)\omega(y,z).\]
By Theorem 3.3. of \cite{Laca_multiplier}, it follows that $\omega$ extends to a multiplier on $\mathbb{R}^{d}$. We denote the extension again by $\omega$. For $x \in \mathbb{R}^{d}$, let $U_{x}$ be the unitary on $L^{2}(\mathbb{R}^{d})$ defined by the following formula
\[
U_{x}f(y)=\omega(x,y-x)f(y-x)\]
for $f \in L^{2}(\mathbb{R}^{d})$. Note that $U_{x}U_{y}=\omega(x,y)U_{x+y}$ for $x,y \in \mathbb{R}^{d}$. Also observe that  for $f,g \in \mathbb{R}^{d}$, the map $\mathbb{R}^{d} \ni x \to \langle U_{x}f|g \rangle$ is measurable. For a proof of this fact, we refer the reader to the paragraph preceding Theorem 3.2. of \cite{Laca_multiplier}.

For $x \in P$, let $\alpha_{x}$ be the automorphism of $B(L^{2}(\mathbb{R}^{d}))$ defined by the formula \[\alpha_{x}(A)=U_{x}AU_{x}^{*}.\] It is clear that $\alpha:=\{\alpha_{x}\}_{x \in P}$ is a semigroup of unital normal $*$-endomorphisms of $B(L^{2}(\mathbb{R}^{d}))$. The weak measurability of $\{U_{x}\}_{x \in P}$ implies that for every $A \in B(L^{2}(\mathbb{R}^{d}))$ and $f,g \in L^{2}(\mathbb{R}^{d})$, the map $P \ni x \to \langle \alpha_{x}(A)f|g \rangle$ is measurable. Now Corollary 4.3. of \cite{Murugan_Sundar} implies that $\alpha$ is an $E_0$-semigroup. 

Let $\displaystyle F:=\coprod_{x \in \Omega}F(x)$ be the product system associated to the $E_0$-semigroup $\alpha$. Then it is clear that for every $x \in \Omega$, $F(x)$ is spanned by $U_{x}$. For $x \in \Omega$, let $\theta_{x}:E(x) \to F(x)$ be the unitary such that $\theta_{x}(e(x))=U_{x}$. Then the map $\displaystyle \theta:=\coprod_{x \in \Omega}\theta_{x}:E \to F$ is $1$-$1$, onto and preserves the multiplication. To see that $\theta$ is measurable, let $\mu:\Omega \times \mathbb{C} \to E$ be defined by $\mu(x,\lambda)=\lambda e(x)$ and let $\nu:\Omega \times \mathbb{C} \to F$ be defined by $\nu(x,\lambda)=(x,\lambda U_{x})$. Then $\mu$ and $\nu$ are measurable. Moreover, $\mu$ is $1$-$1$ and onto. Since the spaces involved are standard, it follows that $\mu^{-1}$ is measurable. Note  that $\theta=\nu \circ \mu^{-1}$. Hence $\theta$ is measurable. This completes the proof. \hfill $\Box$

\begin{rmrk}
\label{endomorphisms and intertwiners}
Let $\clh$ be an infinite dimensional separable Hilbert space. The following is well known. See for instance, Section 2.2 of \cite{Izumi}.  
\begin{enumerate}
\item[(1)] Let $\alpha$ be a normal $*$-endomorphism of $B(\clh)$. The intertwining space of $\alpha$ denoted $\mathcal{E}_{\alpha}$ is defined as 
 \[
 \mathcal{E}_{\alpha}:=\{T \in B(\clh): \alpha(A)T=TA ~\textrm{for~}A \in B(\clh)\}.\]
 Note that if $T,S \in \mathcal{E}_{\alpha}$ then $T^{*}S$ is a scalar, for it commutes with every element of $B(\clh)$, which we denote by $\langle S|T \rangle$. Then $\langle~|~\rangle$ is an inner product on $\mathcal{E}_{\alpha}$. The norm on $\mathcal{E}_{\alpha}$ induced by the inner product and the operator norm coincide. Hence it follows that $\mathcal{E}_{\alpha}$ is a Hilbert space. Moreover $\mathcal{E}_{\alpha}$ is separable. Let $\beta$ be another normal $*$-endomorphism of $B(\clh)$. Then  $\alpha=\beta$ if and only if $\mathcal{E}_{\alpha}=\mathcal{E}_{\beta}$. 
 
  Let $\alpha$ and $\beta$ be normal $*$-endomorphisms of $B(\clh)$. Then $\mathcal{E}_{\alpha \circ \beta}$ is the closed linear span of $\{ST: S \in \mathcal{E}_{\alpha}, T \in \mathcal{E}_{\beta}\}$. 

 \item[(2)] Conversely, let $E \subset B(\clh)$ be a separable norm closed subspace of $B(\clh)$ such that $T^{*}S$ is a scalar for every $S,T \in E$.  For $S,T \in E$, let $\langle S|T \rangle=T^{*}S$. Note that $E$ is a separable Hilbert space w.r.t. to the inner product $\langle~|~\rangle$. Then there exists a unique normal $*$-endomorphism $\alpha$ of $B(\clh)$ such that $E=\mathcal{E}_{\alpha}$. The endomorphism $\alpha$ has the following expression. Let $\{V_{i}\}_{i=1}^{d}$ be an orthonormal basis for $E$ where $d$ denotes the dimension of $E$. Then $\alpha$ is given by the equation 
 \begin{equation}
 \label{definition of alpha}
 \alpha(A)=\sum_{i=1}^{d}V_{i}AV_{i}^{*}.
 \end{equation}
When $d$ is infinite, the sum in Equation \ref{definition of alpha} is a strongly convergent sum.  Moreover  $\alpha(1)$ is the projection onto the closure of $E\clh:=span\{T\xi: T \in E, \xi \in \clh\}$. The endomorphism $\alpha$ is unital if and only if $\sum_{i=1}^{d}V_{i}V_{i}^{*}=1$ if and only if $\overline{E\clh}=\clh$.  
\end{enumerate}
\end{rmrk}

From now on, we assume that the fibres of an abstract product system are infinite dimensional.  Let $\clh$ be an infinite dimensional separable Hilbert space. We endow $B(\clh)$ with the measurable structure induced by the weak operator topology.
\begin{dfn}
\label{representation}
Let $\displaystyle E:=\coprod_{x \in \Omega}E(x)$ be a product system over $\Omega$. By a representation of $E$ on $\clh$, we mean a map 
$\phi:E \to B(\clh)$ such that 
\begin{enumerate}
\item[(1)] the map $\phi$ is measurable, 
\item[(2)] for $u,v \in E$, $\phi(uv)=\phi(u)\phi(v)$, and
\item[(3)] for $x \in \Omega$ and $u,v \in E(x)$, $\phi(v)^{*}\phi(u)=\langle u|v \rangle$.
\end{enumerate}
The representation $\phi$ is called essential if $\overline{\phi(E(x))\clh}=\clh$ for every $x \in \Omega$. 
\end{dfn}
Let $\phi:E \to B(\clh)$ be a representation. Then $\phi$ restricted to each fibre is linear. The proof is exactly the same as in the $1$-dimensional case and hence we omit the proof. For the proof in the $1$-dimensional case, we refer the reader to Page 71 of \cite{Arveson}. Moreover Condition $(3)$ implies that $\phi$ restricted to each fibre is isometric.

Fix $x \in \Omega$. Note that $\phi(E(x))$ is a separable, norm closed subspace of $B(\clh)$ such that $T^{*}S$ is a scalar for every $S,T \in \phi(E(x))$. Hence by  Remark \ref{endomorphisms and intertwiners}, it follows that there exists a unique normal $*$-endomorphism $\alpha_{x}$ of $B(\clh)$ such that \[
\phi(E(x))=\{T \in B(\clh): \alpha_{x}(A)T=TA ~\textrm{for~}A \in B(\clh)\}.\]
Recall from Remark \ref{endomorphisms and intertwiners} that $\alpha_{x}(1)$ is the projection onto the closed subspace $\overline{\phi(E(x))\clh}$.

Let $x,y \in \Omega$ be given. Since $E(x+y)$ is the closure of the  linear span of the set $\{uv:u \in E(x),v \in E(y)\}$, it follows from $(2)$ of Definition \ref{representation} that $\phi(E(x+y))$ is the closed linear span of $\{ST: S \in \phi(E(x)), T \in \phi(E(y))\}$. Hence by Remark \ref{endomorphisms and intertwiners}, it follows that $\alpha_{x+y}=\alpha_{x} \circ \alpha_{y}$. Consequently $\{\alpha_{x}\}_{x \in \Omega}$ is a semigroup of normal $*$-endomorphisms of $B(\clh)$. 
 Now suppose that $\phi$ is essential. Then by Remark \ref{endomorphisms and intertwiners}, it follows that $\alpha_{x}$ is unital for every $x \in \Omega$. 

Since $E \simeq \Omega \times \ell^{2}$, it follows that there exists measurable sections $e_1,e_2,\cdots$ such that for every $x \in \Omega$, $\{e_{1}(x),e_{2}(x),\cdots\}$ forms an orthonormal basis for $E(x)$. Consequently, for $x \in \Omega$, $\{\phi(e_1(x)),\phi(e_2(x)),\cdots\}$ is an orthonormal basis for $\phi(E(x))$. Hence by Remark \ref{endomorphisms and intertwiners}, it follows that for $x \in \Omega$, $\alpha_{x}$ is given by the equation
\begin{equation}
\label{defn of alpha}
\alpha_{x}(A)=\sum_{n=1}^{\infty}\phi(e_n(x))A\phi(e_n(x))^{*}
\end{equation}
where the sum in Equation \ref{defn of alpha} is a strongly convergent sum. The measurability of $\phi$ and Equation \ref{defn of alpha} implies that for $A \in B(\clh)$ and $\xi,\eta \in \clh$, the map $\Omega \ni x \to \langle \alpha_{x}(A)\xi|\eta \rangle \in \mathbb{C}$ is measurable. By Proposition 4.2 of \cite{Murugan_Sundar}, it follows that for $A \in B(\clh)$ and $\xi,\eta \in \clh$, the map $\Omega \ni x \to \langle \alpha_{x}(A)\xi|\eta \rangle \in \mathbb{C}$ is continuous. Again by Proposition 4.2 of \cite{Murugan_Sundar}, it follows that $\{\alpha_{x}\}_{x \in \Omega}$ extends to a unique $E_0$-semigroup which we still denote by $\alpha:=\{\alpha_{x}\}_{x \in P}$. The constructed $E_0$-semigroup $\alpha$ is called the $E_0$-semigroup associated to the essential representation $\phi$. 

\begin{ppsn}
\label{Essential}
Let $\phi:E \to B(\clh)$ be an essential representation and let $\alpha:=\{\alpha_{x}\}_{x \in P}$ be the $E_0$-semigroup associated to $\phi$. Then $E$ is isomorphic to the product system associated to $\alpha$.
\end{ppsn}
\textit{Proof.} Let $\displaystyle F:=\coprod_{x \in \Omega}F(x)$ be the product system associated to $\alpha$. For $x \in \Omega$, by the definition of $\alpha_{x}$, $F(x)=\phi(E(x))$. Now the map $E \ni u \to (p(u),\phi(u)) \in F$ is an isomorphism of product systems. Here $p:E \to \Omega$ is the canonical surjection that comes equipped with the product system $E$. This completes the proof. \hfill $\Box$

\section{Construction of an essential representation}
We fix an element $a \in \Omega$ for the rest of this section. For $x,y \in \mathbb{R}^{d}$, we write $x>y$ if $x-y \in \Omega$. We have the following archimedean principle. 
\begin{lmma}
\label{archimedean}
Let $x \in \mathbb{R}^{d}$ be given. Then there exists a positive integer $n_0$ such that $n_0a >x$.
\end{lmma}
\textit{Proof.} Note that the sequence $a-\frac{x}{n} \to a \in \Omega$. But $\Omega$ is an open subset of $\mathbb{R}^{d}$.  Hence there exists a positive integer $n_0$ such that $n \geq n_0$ implies $a-\frac{x}{n} \in \Omega$. This implies that $na-x \in \Omega$ for $n \geq n_0$. In particular, $n_0a-x \in \Omega$. This completes the proof.  \hfill $\Box$

\begin{lmma}
\label{disjoint}
The intersection $\displaystyle \bigcap_{n=0}^{\infty}(\Omega+na)=\emptyset$. Also $\{\Omega+na\}_{n=0}^{\infty}$ is a decreasing sequence of subsets of $\Omega$. 
\end{lmma}
\textit{Proof.} Suppose $y \in \displaystyle \bigcap_{n=0}^{\infty}(\Omega+na)$. Then $y-na \in \Omega$ for every $n \geq 0$. By Lemma \ref{archimedean}, it follows that there exists a positive integer $n_0$ such that $n_0a-y \in \Omega$. Now observe that $-a=(y-(n_0+1)a)+(n_0a-y) \in \Omega$ which is a contradiction since $\Omega \cap -\Omega=\emptyset$. It is clear that $\{\Omega+na\}_{n=0}^{\infty}$ is a decreasing sequence of subsets of $\Omega$.  \hfill $\Box$

Let us fix notation that we will use throughout  this paper.   Let \[L_{k}:=(\Omega+ka)\backslash(\Omega+(k+1)a)\] for $k \in \mathbb{N}$. Then Lemma \ref{disjoint} implies that $\{L_{k}:k \in \mathbb{N}\}$ is a disjoint family of measurable subsets of $\Omega$. Note that for $k \in \mathbb{N}$, $\displaystyle \Omega+ka=\coprod_{m \geq k}L_{m}$. Since $\displaystyle \Omega=\coprod_{k \in \mathbb{N}}L_{k}$, it follows that given $x \in \Omega$, there exists a unique non-negative integer $n(x)$ such that $x \in L_{n(x)}$. Since $n(x)=k$ for $x \in L_{k}$, it is clear that the map $\Omega \ni x \to n(x) \in \mathbb{N}$ is measurable.  Note that $n(x+a)=n(x)+1$ for $x \in \Omega$.  Also observe that for $x \in \Omega$, $x-n(x)a \in L_0$ and for $x \in L_0$ and $k \in \mathbb{N}$, $x+ka \in L_{k}$. 

We need the fact that $L_{k}$ has non-zero Lebesgue measure for every $k \in \mathbb{N}$. In what follows, we will not use any other measure except the Lebesgue measure on $\mathbb{R}^{d}$. Also we denote the Lebesgue measure on $\mathbb{R}^{d}$ by $\lambda$. Since $L_{k}=L_0+ka$, it suffices to show that $L_0=\Omega\backslash(\Omega+a)$ has positive measure. Observe that $\Omega \backslash(\Omega+a)$ contains the open set $\Omega\backslash (P+a)$ which is non-empty since $\frac{a}{2} \in \Omega \backslash (P+a)$. Now it follows immediately that $\Omega\backslash(\Omega+a)$ has positive measure. For $z \in \mathbb{R}^{d}$, let $L_{z}:=(L_0+z)\cap \Omega$.

Let $p:E \to \Omega$ be a product system which is fixed for the rest of this section. We assume that the fibres of $E$ are infinite dimensional. Let $e \in E(a)$ be a unit vector which is fixed for the rest of this section. Our goal is to exhibit an essential representation of $E$ on an infinite dimensional separable Hilbert space. 

Let $\mathcal{V}$ denote the vector space of measurable sections of $E$ which are square integrable over $L_{z}$ for every $z \in \mathbb{R}^{d}$. More precisely, let $f:\Omega \to E$ be a measurable section. Then $f \in \mathcal{V}$ if and only if for every $z \in \mathbb{R}^{d}$, $\displaystyle \int_{L_{z}}||f(x)||^{2}dx < \infty$. 

Let $f \in \mathcal{V}$ and $k \in \mathbb{N}$ be given. We say that $f$ is $k$-stable if $f(x+a)=f(x)e$ for almost all $x>ka$, i.e. the measurable set $\{x \in \Omega+ka: f(x+a) \neq f(x)e\}$ has measure zero. Note that if $f$ is $k$-stable and $k_1 \geq k$ then $f$ is $k_1$-stable. We say a section in $\mathcal{V}$ is stable if it is $k$-stable for some $k \in \mathbb{N}$. Denote the set of stable sections in $\mathcal{V}$ by $\mathcal{S}$. Then it is clear that $\mathcal{S}$ is a vector subspace of $\mathcal{V}$. 

Let $f \in \mathcal{V}$. We say that $f$ is null if there exists $k \in \mathbb{N}$ such that $f(x)=0$ for almost all $x>ka$, i.e. there exists $k \in \mathbb{N}$ such that $\{x>ka: f(x) \neq 0\}$ has measure zero. Denote the set of null sections in $\mathcal{V}$ by $\mathcal{N}$. Then it is clear that $\mathcal{N}$ is a vector subspace of $\mathcal{V}$. We leave it to the reader to verify that $\mathcal{N} \subset \mathcal{S}$.

\begin{lmma}
\label{m stability}
Let $f \in \mathcal{S}$ be given. Assume that $f$ is $k$-stable for some $k \in \mathbb{N}$. Then for every $m \geq 1$, $f(x+ma)=f(x)e^{m}$ for almost all $x>ka$. 
\end{lmma}
\textit{Proof.} We prove this by induction on $m$. Let $A_{m}:=\{x>ka: f(x+ma)\neq f(x)e^{m}\}$. The fact that $f$ is $k$-stable implies that $A_1$ has measure zero. Now assume that $A_m$ has measure zero. Let $x>ka$ be given. Suppose $x \notin A_m$ and $x+ma \notin A_1$. Then calculate as follows to observe that 
\begin{align*}
f(x+(m+1)a)&=f(x+ma+a) \\
                   &=f(x+ma)e ~~(\textrm{Since $x+ma \notin A_1$}) \\
                   &=f(x)e^{m}e ~~(\textrm{Since $x \notin A_m$}) \\
                   &=f(x)e^{m+1}.
\end{align*}
This implies that $A_{m+1} \subset A_{m} \cup ( (\Omega+ka)\cap (A_1-ma)) \subset A_m \cup (A_1-ma)$. Since $A_m$ and $A_1$ have measure zero, it follows that $A_{m+1}$ has measure zero.  This completes the proof. \hfill $\Box$

Let $f,g \in \mathcal{S}$ be given. Since $f$ and $g$ are square integrable over $L_{z}$ for every $z \in \mathbb{R}^{d}$, it follows that the integral $\displaystyle \int_{L_k}\langle f(x)|g(x)\rangle dx$ exists for every $k \in \mathbb{N}$.

\begin{lmma}
\label{inner product1}
Let $f,g \in \mathcal{S}$ be given. Assume that $f$ and $g$ are $k_0$-stable for some $k_0 \in \mathbb{N}$. Then for $k \geq k_0$, 
\[
\int_{L_k}\langle f(x)|g(x)\rangle dx \rangle=\int_{L_{k_0}}\langle f(x)|g(x) \rangle dx.\]
\end{lmma}
\textit{Proof.} Let $k>k_0$ be given. Note that the map $L_{k_0} \ni x \to x+(k-k_0)a \in L_{k}$ is a measurable bijection which is measure preserving. Calculate as follows to observe that 
\begin{align*}
\int_{L_k}\langle f(x)|g(x) \rangle dx&=\int_{L_{k_0}}\langle f(x+(k-k_0)a)|g(x+(k-k_0)a)\rangle dx \\
                                                        & = \int_{L_{k_0}}\langle f(x)e^{k-k_0}|g(x)e^{k-k_0}\rangle dx ~~(\textrm{by Lemma \ref{m stability}})\\
                                                        &= \int_{L_{k_0}}\langle f(x)|g(x) \rangle .
                                                        \end{align*}
This completes the proof. \hfill $\Box$

Let $f,g \in \mathcal{S}$ be given. Thanks to Lemma \ref{inner product1},  $\displaystyle \lim_{k \to \infty} \int_{L_k}\langle f(x)|g(x)\rangle dx$ exists. Define \[
\langle f|g \rangle = \lim_{k \to \infty} \int_{L_k}\langle f(x)|g(x)\rangle dx.\]  Observe that $\langle~|~\rangle$ defines a semi-definite inner product on $\mathcal{S}$. Let $f \in \mathcal{S}$. We claim that $\langle f|f \rangle=0$ if and only if $f \in \mathcal{N}$. It is clear that if $f \in \mathcal{N}$, then $\langle f|f\rangle=0$. Now let $f \in \mathcal{S}$ be such that  $\langle f|f\rangle=0$.  Assume that $f$ is $k_0$-stable for some $k_0 \in \mathbb{N}$. Lemma \ref{inner product1} implies that $\displaystyle \int_{L_k}||f(x)||^{2}dx=0$ for every $k \geq k_0$. Consequently for every $k \geq k_0$, $f(x)=0$ for almost all $x \in L_{k}$. Since $\Omega+k_0a=\displaystyle \coprod_{k \geq k_0}L_{k}$, it follows that $f(x)=0$ for almost all $x>k_0a$. This proves that $f \in \mathcal{N}$.  

Thus $\langle~|~\rangle$ descends to a positive definite inner product on $\mathcal{S}/\mathcal{N}$ which we still denote by $\langle~|~\rangle$. Let $\clh$ be the completion of the pre-Hilbert space $\mathcal{S}/\mathcal{N}$.


\begin{ppsn}
The Hilbert space $\clh$ is separable and is infinite dimensional.
\end{ppsn}
\textit{Proof.} Let $k \in \mathbb{N}$ be fixed. Consider a measurable section $\xi:L_{k} \to E$ which is square integrable. Define a section $\widetilde{\xi}:\Omega \to E$ by the following formula
\begin{equation*}
 \widetilde{\xi}(x):=\begin{cases}
 \xi(x-n(x)a+ka)e^{n(x)-k} & \mbox{ if
} x > ka ,\cr
    & \cr
    0 & \mbox{elsewhere}.
         \end{cases}
\end{equation*}
Note that for $x \in L_{k}$, $\widetilde{\xi}(x)=\xi(x)$. For $m > k$ and   $x \in L_{m}$, $\widetilde{\xi}(x)=\xi(x-ma+ka)e^{m-k}$. This implies that $\widetilde{\xi}$ is measurable on each $L_m$ for $m \geq k$. Hence $\widetilde{\xi}$ is measurable on $\displaystyle \coprod_{m \geq k} L_{m}=\Omega+ka$. It is clear that $\widetilde{\xi}$ is measurable on the complement of $\Omega+ka$. Consequently it follows that $\widetilde{\xi}$ is a measurable section.

We claim that $\widetilde{\xi} \in \mathcal{V}$.  Let $z \in \mathbb{R}^{d}$ be given. Set $A:=L_{z} \cap (\Omega+ka)$. Since $\widetilde{\xi}$ vanishes on $L_{z}\backslash A$, it follows that  $\displaystyle \int_{L_{z}}||\widetilde{\xi}(x)||^{2}dx=\int_{A}||\widetilde{\xi}(x)||^{2}dx$. If $A$ is empty, then the last equality implies that $\displaystyle \int_{L_{z}}||\widetilde{\xi}(x)||^{2}dx=0 <\infty$. In other words, $\widetilde{\xi}$ is square integrable on $L_{z}$. Suppose that $A$ is non-empty.  Let $\chi:A \to L_{k}$ be the map defined by \[\chi(x)=x-n(x)a+ka.\] The measurability of the map $\Omega \ni x \to n(x) \in \mathbb{N}$ implies that $\chi$ is measurable. We claim that $\chi$ is $1$-$1$. 
Let $x_1,x_2 \in A$ be such that $\chi(x_1)=\chi(x_2)$. To prove $x_1=x_2$, it suffices to show that $n(x_1)=n(x_2)$. Suppose not. Without loss of generality, we can assume that $n(x_1)<n(x_2)$.  Note that $x_1-z \in \Omega$ and $n(x_2)-n(x_1)$ is a natural number greater than or equal to $1$. Since $\{\Omega+na: n \geq 1\}$ is a decreasing sequence of subsets, it follows that $x_2-z=(x_1-z)+(n(x_2)-n(x_1))a \in \Omega+a$ which is a contradiction to the fact that $x_2-z \in L_0 = \Omega\backslash(\Omega+a)$. This contradiction implies that $n(x_1)=n(x_2)$ and consequently $x_1=x_2$. This proves that $\chi$ is $1$-$1$. Since $A$ and $L_{k}$ are $G_{\delta}$ subsets of $\mathbb{R}^{d}$, it follows that $A$ and $L_{k}$ are Polish spaces. Let $B$ be the image of $\chi$. Then by Theorem 3.3.2. of \cite{Arveson_invitation}, it follows that $B$ is a Borel subset of $L_{k}$ and $\chi:A \to B$ is a Borel isomorphism. 

We claim that $\chi$ is measure preserving. Let $C \subset A$ be a Borel subset. For $m \geq k$, let $C_{m}:=\{x \in C: n(x)=m\}$. Then $\displaystyle C=\coprod_{m \geq k}C_{m}$. As a consequence, we have $\displaystyle \chi(C)=\coprod_{m \geq k}\chi(C_{m})=\coprod_{m \geq k}(C_{m}-ma+ka)$. Now calculate as follows to observe that 
\begin{align*}
\lambda(\chi(C))&=\sum_{m \geq k} \lambda(C_{m}-ma+ka) \\
                          &=\sum_{m \geq k}\lambda(C_{m}) \\
                          &= \lambda(\coprod_{m \geq k}C_{m}) \\
                          &=\lambda(C).
\end{align*}
This shows that $\chi$ is measure preserving. Calculate as follows to observe that 
\begin{align*}
\int_{L_{z}}||\widetilde{\xi}(x)||^{2}dx&= \int_{A}||\widetilde{\xi}(x)||^{2}dx\\
                                                     &= \int_{A}||\xi(x-n(x)a+ka)||^{2} dx \\
                                                   & =\int_{A}||\xi(\chi(x))||^{2}dx \\
                                                   &=\int_{B}||\xi(x)||^{2}dx \\
                                                   &\leq \int_{L_k}||\xi(x)||^{2}dx\\
                                                   &< \infty.
\end{align*}
This shows that $\widetilde{\xi} \in \mathcal{V}$. Next we claim $\widetilde{\xi}$ is $k$-stable. Let $x>ka$ be given. Calculate as follows to observe that 
\begin{align*}
\widetilde{\xi}(x+a)&=\xi(x+a-n(x+a)a+ka)e^{n(x+a)-k} \\
                             &=\xi(x+a-(n(x)+1)a+ka)e^{n(x)+1-k} ~~(\textrm{~Since $n(x+a)=n(x)+1$})\\
                             &=\xi(x-n(x)a+ka)e^{n(x)-k}e \\
                             &=\widetilde{\xi}(x)e.
\end{align*}
This proves that $\widetilde{\xi}$ is $k$-stable. Note that $\widetilde{\xi}(x)=\xi(x)$ for $x \in L_{k}$. Let $\clh_{k}=L^{2}(L_k,E)$. By Lemma \ref{inner product1}, it follows that the map $\clh_{k} \ni \xi \to \widetilde{\xi}+\mathcal{N}$ is well-defined and is an isometry which we denote  by $V_{k}$.

Let $f \in \mathcal{V}$ be given. Assume that $f$ is $k$-stable for some $k \in \mathbb{N}$. Let $\xi:L_{k} \to E$ be defined by $\xi(x)=f(x)$ for $x \in L_{k}$. For $m \geq 1$, let $A_{m}:=\{x>ka: f(x+ma)\neq f(x)e^{m}\}$.  By Lemma \ref{m stability}, it follows that $A_{m}$ has measure zero for every $m \geq 1$. Let $m>k$ and $x \in L_{m}$ be given. Suppose $x \notin A_{m-k}+(m-k)a$. Then calculate as follows to observe that 
\begin{align*}
f(x)&=f(x-ma+ka+(m-k)a) \\
     &= f(x-ma+ka)e^{m-k} ~~(\textrm{Since $x-ma+ka \notin A_{m-k}$ and $x-ma+ka >ka$.})\\ 
      &=\widetilde{\xi}(x).
\end{align*}
Therefore  $\{x \in L_{m}: f(x) \neq \widetilde{\xi}(x)\}\subset A_{m-k}+(m-k)a$. Since $A_{m-k}+(m-k)a$ has measure zero, it follows that $\{x \in L_{m}: f(x) \neq \widetilde{\xi}(x)\}$ has measure zero.  Thus for every $m>k$, $f(x)=\widetilde{\xi}(x)$ for almost all $x \in L_{m}$. By definition, $f$ agrees with $\widetilde{\xi}$ on $L_{k}$.  Hence $f(x)=\widetilde{\xi}(x)$ for almost all $x \in \displaystyle \coprod_{m \geq k}L_{m}=\Omega+ka$. Consequently $f-\widetilde{\xi} \in \mathcal{N}$. This proves that $f+\mathcal{N}=\widetilde{\xi}+\mathcal{N}$. Thus we have shown that $\{f+\mathcal{N}: f \in \mathcal{S}\}= \bigcup_{k=0}^{\infty}V_{k}\clh_{k}$. Since each $\clh_k$ is separable, it follows that $\clh$ is separable. As each $\clh_{k}$ is infinite dimensional, it follows that $\clh$ is infinite dimensional. This completes the proof. \hfill $\Box$

We need the following two important lemmas before defining a representation of $E$ on $\clh$. Fix $k \in \mathbb{N}$. Let $b \in \Omega$ be such that $b>ka$. Recall that $L_{k}=(\Omega+ka)\backslash(\Omega+(k+1)a)$ and $L_{b}=(\Omega+b)\backslash(\Omega+b+a)$. Let $x \in L_{k}$ be given. By Lemma \ref{archimedean}, there exists $m_0 \in \mathbb{N}$ such that $m_{0}a-(b-x)=x+m_{0}a-b \in \Omega$ or in other words $x+m_0a \in \Omega+b$.  Let $m(x)$ be the least non-negative integer such that $x+m(x)a \in \Omega+b$. 

\begin{lmma}
\label{bijection}
With the foregoing notations, we have the following. 
\begin{enumerate}
\item[(1)] For every $x \in L_{k}$, $x+m(x)a \in L_{b}$.
\item[(2)] For every $x \in L_{k}$, the intersection $\{x+ma: m \in \mathbb{N}\} \cap L_{b}$ is singleton.
\item[(3)] The map $\chi:L_{k} \to L_{b}$ defined by $\chi(x)=x+m(x)a$ is a measurable bijection.
\item[(4)] The map $\chi:L_{k} \to L_{b}$ is measure preserving. 
\end{enumerate}
\end{lmma}
\textit{Proof.} Let $x \in L_{k}$ be given. Suppose $m(x)=0$. Then since $x\ngtr(k+1)a$ and $b>ka$, it follows that $x \ngtr b+a$. Hence $x=x+m(x)a \in L_{b}$. Suppose $m(x) \geq 1$. Then by definition $x+(m(x)-1)a \notin \Omega+b$ or in other words, $x+m(x)a \notin \Omega+b+a$. In this case too, $x+m(x)a \in L_{b}$. This proves $(1)$. 

Let $x \in L_{k}$ be given. Suppose $x+ma \in L_{b}$ for some $m \in \mathbb{N}$. Then $x+ma \in \Omega+b$. Hence by the definition of $m(x)$, it follows that $m \geq m(x)$. To prove $(2)$, it suffices to show that $m=m(x)$. Suppose not. Then $m >m(x)$. Write $m=m(x)+n$ with $n \geq 1$. Now $x+ma=(x+m(x)a)+na \in (\Omega+b)+na \subset \Omega+b+a$ which is a contradiction since $x+ma \notin \Omega+b+a$. This contradiction proves that $m=m(x)$. This proves $(2)$. 

Let $x_1,x_2 \in L_{k}$ be such that $\chi(x_1)=\chi(x_2)$. Then $x_1+m(x_1)a=x_2+m(x_2)a$. To show $x_1=x_2$, it suffices to show that $m(x_1)=m(x_2)$. Suppose not. Without loss of generality, we can assume that $m(x_1)<m(x_2)$. Then $x_{1}=x_{2}+(m(x_2)-m(x_1))a \in \Omega+ka+(m(x_2)-m(x_1))a \subset \Omega+(k+1)a$ which is a contradiction since $x_1 \ngtr (k+1)a$. This contradiction implies that $m(x_1)=m(x_2)$ and consequently $x_1=x_2$. This proves that $\chi$ is $1$-$1$. 

Let $y \in L_{b}$ be given. Then $y>b>ka$. Hence the set $\{n \in \mathbb{N}: y-na \in \Omega+ka\}$ is non-empty, for it contains $0$. We claim that the set $\{n \in \mathbb{N}: y-na \in \Omega+ka\}$ is bounded. Suppose not. Then there exists a sequence $(n_\ell)$ such that $n_{\ell} \to \infty$ and $y-n_{\ell}a \in \Omega+ka \subset P$. Hence $\frac{y}{n_\ell}-a \in P$ for every $\ell$. But $\frac{y}{n_\ell}-a \to -a$. This forces that $-a \in P$ or $a \in -P$ which is a contradiction.  Let $n_0$ be the largest non-negative integer such that $y-n_0a \in \Omega+ka$. Then $y-(n_0+1)a \notin \Omega+ka$ or in other words $y-n_0a \notin \Omega+(k+1)a$.  Let $x:=y-n_0a$. Then $x \in L_{k}$ and $y=x+n_0a \in L_{b}$. Since the intersection $\{x+ma: m \in \mathbb{N}\}\cap L_{b}$ is singleton, it follows that $y=\chi(x)$. This proves that $\chi$ is onto. 

To show $\chi$ is measurable, it is enough to show that $L_{k} \ni x \to m(x) \in \mathbb{N}$ is measurable. Let $r  \in \mathbb{R}$ be given. We claim that $\{x \in L_{k}: m(x) \geq r \}$ is a closed subset of $L_{k}$. Let $(x_n)$ be a sequence in $L_{k}$ such that $m(x_n) \geq r$ and  $x_n \to x \in L_{k}$. Then the sequence $x_{n}+m(x)a \to x+m(x)a \in \Omega+b$. But $\Omega+b$ is an open subset of $\mathbb{R}^{d}$ containing $x+m(x)a$. Hence $x_n+m(x)a \in \Omega+b$ eventually. By the definition of the function $m$, it follows that $m(x_n) \leq m(x)$ eventually. Thus $m(x) \geq r$. This proves that $\{x \in L_{k}: m(x) \geq r\}$ is a closed subset of $L_{k}$. As a consequence, we obtain that the function $m$ is measurable and consequently $\chi$ is measurable. This proves $(3)$. 

Since $L_{k}$ and $L_{b}$ are $G_{\delta}$-subsets of $\mathbb{R}^{d}$, it follows that $L_{k}$ and $L_{b}$ are Polish spaces. Hence by Theorem 3.3.2. of \cite{Arveson_invitation}, it follows that $\chi$ is a Borel isomorphism. Let $A \subset L_{k}$ be a measurable subset. For $n \in \mathbb{N}$, let $A_{n}:=\{x \in A: m(x)=n\}$. Then $\displaystyle A=\coprod_{n=0}^{\infty}A_{n}$ and $\displaystyle \chi(A)=\coprod_{n=0}^{\infty}\chi(A_n)=\coprod_{n=0}^{\infty}(A_n+na)$. Now calculate as follows to observe that 
\begin{align*}
\lambda(\chi(A))&=\sum_{n=0}^{\infty}\lambda(A_n+na) \\
                         &=\sum_{n=0}^{\infty}\lambda(A_n) \\
                         &=\lambda(\coprod_{n=0}^{\infty}A_{n})\\
                         &=\lambda(A).
\end{align*}
This proves $(4)$. This completes the proof. \hfill $\Box$

\begin{lmma}
\label{shift}
Let $f,g \in \mathcal{S}$ be given. Assume that $f$ and $g$ are $k$-stable for some $k \in \mathbb{N}$. Let $b \in \Omega$ be such that $b>ka$. Then 
\[
\langle f|g \rangle = \int_{L_b}\langle f(x)|g(x) \rangle dx.\]
\end{lmma}
\textit{Proof.} Let $L_{k} \ni x\to m(x) \in \mathbb{N}$ and $\chi:L_{k} \to L_{b}$ be the functions considered in Lemma \ref{bijection}. For $n \in \mathbb{N}$, let $A_{n}:=\{x \in L_{k}:m(x)=n\}$. Then $L_{k}=\displaystyle \coprod_{n=0}^{\infty}A_{n}$.
Now calculate as follows to observe that 
\begin{align*}
\int_{L_{b}}\langle f(x)|g(x) \rangle dx & = \int_{L_{k}}\langle f(\chi(x))|g(\chi(x))\rangle dx ~~(\textrm{Since $\chi$ is measure preserving}) \\
                                                           & = \sum_{n=0}^{\infty}\int_{A_n}\langle f(\chi(x))|g(\chi(x))\rangle dx \\
                                                           &= \sum_{n=0}^{\infty}\int_{A_n} \langle f(x+na)|g(x+na) \rangle dx \\
                                                           &=\sum_{n=0}^{\infty}\int_{A_n}\langle f(x)e^{n}|g(x)e^{n} \rangle dx ~~(\textrm{by Lemma \ref{m stability}})\\
                                                           &=\sum_{n=0}^{\infty}\int_{A_n}\langle f(x)|g(x) \rangle dx \\
                                                           &=\int_{L_{k}}\langle f(x)|g(x) \rangle dx \\
                                                           &=\langle f|g \rangle ~~(\textrm{by Lemma \ref{inner product1}}) .
                                                           \end{align*}
This completes the proof. \hfill $\Box$

Let $b \in \Omega$ and $v \in E(b)$ be given. For  $f \in \mathcal{S}$, let $\phi_{0}(v)f:\Omega \to E$ be the measurable section defined by 
\begin{equation*}
 (\phi_{0}(v)f)(x):=\begin{cases}
 vf(x-b) & \mbox{ if
} x > b \cr
    & \cr
    0 & \mbox{elsewhere}.
         \end{cases}
\end{equation*}
Let $f \in \mathcal{S}$ be given. We leave it to the reader to verify that $\phi_{0}(v)f \in \mathcal{V}$. Assume that $f$ is $k$-stable for some $k \geq 1$. Set $A:=\{x>ka: f(x+a) \neq f(x)e\}$. Then $A$ has measure zero. Choose $k_0 \in \mathbb{N}$ such that $k_0a>b$ and set $k_1=k_0+k$. We claim that $\phi_{0}(v)f$ is $k_1$-stable.  Let $x>k_1a$ and $x \notin A+b$ be given. Calculate as follows to observe that 
\begin{align*}
(\phi_{0}(v)f)(x+a)&=vf(x+a-b) ~~(\textrm{Since $x+a >k_1a>k_0a>b$})\\
                            &=vf(x-b)e ~~(\textrm{Since $x-b>k_1a-b>ka$ and $x-b \notin A$}) \\
                            &=(\phi_{0}(v)f)(x)e.
\end{align*}
Hence the set $\{x>k_1a:(\phi_{0}(v)f)(x+a) \neq (\phi_{0}(v)f)(x)e\}$ is contained in $A+b$ which has measure zero. This proves that $\phi_{0}(v)f$ is $k_1$-stable. 

\begin{ppsn}
\label{repn}
Let $b \in \Omega$ and $u,v \in E(b)$ be given. Then for $f \in \mathcal{S}$, \[\langle \phi_{0}(u)f|\phi_{0}(v)f\rangle=\langle u|v \rangle \langle f|f \rangle.\] 
\end{ppsn}
\textit{Proof.} Assume that $f$ is $k$-stable for some $k \geq 1$. Choose $k_0 \geq 1$ such that $k_0a>b$ and set $k_1=k_0+k$. Now calculate as follows to observe that 
\begin{align*}
\langle \phi_{0}(u)f|\phi_{0}(v)f \rangle &= \int _{L_{k_1}}\langle (\phi_{0}(u)f)(x)|(\phi_{0}(v)f)(x) \rangle dx \\
                                                             &= \int_{L_{k_1}}\langle u f(x-b)|v f(x-b) \rangle dx \\
                                                             &=\langle u|v \rangle \int_{L_{k_1}}\langle f(x-b)|f(x-b)\rangle dx \\
                                                             &=\langle u|v \rangle \int_{L_{k_1a-b}}\langle f(x)|f(x) \rangle dx \\
                                                             &= \langle u|v\rangle \langle f|f \rangle ~~(\textrm{since $k_1a-b>ka$ and by Lemma \ref{shift}}).
\end{align*}
This completes the proof. \hfill $\Box$

Let $b \in \Omega$  and $v \in E(b)$ be given.  Proposition \ref{repn} implies that for $f \in \mathcal{S}$, \[\langle  \phi_{0}(v)f|\phi_{0}(v)f \rangle=||v||^{2}\langle  f|f\rangle.\] As a consequence, it follows that there exists a unique bounded linear operator, denoted $\phi(v)$, on $\clh$ such that 
$\phi(v)(f+\mathcal{N})=\phi_{0}(v)f+\mathcal{N}$ for every $f \in \mathcal{S}$. Prop. \ref{repn} implies that for $u,v \in E(b)$, $\phi(v)^{*}\phi(u)= \langle u|v\rangle$. It is clear that $\phi:E \to B(\clh)$ is multiplicative. Next we verify that $\phi$ is measurable. 
We need the following remark, which follows from Tonelli's theorem. 

\begin{rmrk}
\label{measurability}
Let $(X,\mathcal{B}_{X})$ be a measurable space and $(Y,\mathcal{B}_{Y},\mu)$ be a $\sigma$-finite measure space. Let $\mathcal{B}_{X} \otimes \mathcal{B}_{Y}$ be the product $\sigma$-algebra on $X \times Y$. Consider a measurable function $f:X \times Y \to \mathbb{C}$. Suppose that for every $x \in X$, the map $Y \ni y \to f(x,y) \in \mathbb{C}$ is integrable.  Then the function $X \ni x \to \int f(x,y)d\mu(y) \in \mathbb{C}$ is measurable.  
\end{rmrk}

\begin{ppsn}
\label{measurability of phi}
The map $\phi:E \to B(\clh)$ is measurable. Hence $\phi$ is a representation of $E$ on $\clh$.
\end{ppsn}
\textit{Proof.} For $k \geq 1$, let $E_{k}:=\{v \in E: 0<p(v)<ka\}$. Then $E_{k}$ is a measurable subset of $E$ for every $k$ and $\bigcup_{k=1}^{\infty}E_{k}=E$. Thus it suffices to show that $\phi$ restricted to $E_{k}$ is measurable for every $k$. Fix $k \geq 1$. It suffices to show that for $f \in \mathcal{S}$, the map  $E_{k} \ni v \to \langle \phi_{0}(v)f|f \rangle \in \mathbb{C}$ is measurable. Let $f \in \mathcal{S}$ be given. Assume that $f$ is $k_0$-stable for some $k_0 \geq 1$. Then for $v \in E_{k}$, $\phi_{0}(v)f$ is $k_0+k$-stable. Hence by Lemma \ref{inner product1}, it follows that 
\[
\langle \phi_{0}(v)f|f \rangle=\int_{L_{k_0+k}}\langle vf(x-p(v))|f(x) \rangle dx.
\]
The above integral representation together with Remark \ref{measurability} implies that the function $E_{k} \ni v \to \langle \phi_{0}(v)(f)|f \rangle \in \mathbb{C}$ is measurable. This completes the proof. \hfill $\Box$


\begin{rmrk}
We need the following before we proceed further. 
\begin{enumerate}
\item[(1)] Let $x,y \in \Omega$ be such that $x<y$. For $v \in E(x)$ and $w \in E(y)$, there exists a unique element denoted $v^{*}w \in E(y-x)$ such that $\langle v^{*}w|u \rangle=\langle w|vu \rangle$ for $u \in E(y-x)$. 
              Note that for $v \in E(x)$ and $w \in E(y)$, \begin{equation}
                              \label{inequality}
                              ||v^{*}w|| \leq ||v||||w||.
                              \end{equation}
               Let $x,y,z \in \Omega$ be such that $x<y$. Let $v \in E(x)$, $w_1 \in E(y)$ and $w_{2} \in E(z)$ be given. Then \begin{equation}
                                                                                                                                                                                                \label{star operation}
                                                                                                                                                                                                v^{*}(w_1w_2)=(v^{*}w_1)w_2.
                                                                                                                                                                                                \end{equation} The proof of this fact is exactly similar to the proof of Lemma 2.4 of \cite{Murugan_Sundar_discrete}. Hence we omit the proof. 
   \item[(2)]  Let $x, y \in \Omega$ be such that $x<y$. Let $\{v_1,v_2,\cdots,\}$ be an orthonormal basis for $E(x)$. Then for $\xi \in E(y)$, \begin{equation}
   \label{onb}
   \sum_{i=1}^{\infty}||v_{i}^{*}\xi||^{2}=||\xi||^{2}.
   \end{equation}         
   The proof of this fact is exactly similar to the proof of Lemma 2.4 of \cite{Arv}.  Hence we omit the proof.

\end{enumerate}
\end{rmrk}

Let $v \in E(a)$ be given. For $f \in \mathcal{S}$, let $f_{v}:\Omega \to E$ be defined by \[f_{v}(x)=v^{*}f(x+a).\]  Let $f \in \mathcal{S}$ be given. Note that $f_{v}$ is a section. To see that $f_{v}$ is measurable, let $s:\Omega \to E$ be a measurable section. Then the map $\Omega \ni x \to \langle f_{v}(x)|s(x)\rangle=\langle f(x+a)|vs(x)\rangle \in \mathbb{C}$ is measurable. This implies that $f_{v}$ is measurable. We leave it to the reader to verify that $f_{v} \in \mathcal{S}$. We only indicate that to prove $f_{v} \in \mathcal{V}$, one needs to use  Inequality \ref{inequality} and to prove $f_{v}$ is stable one needs to appeal to Equation \ref{star operation}. Note that if $f$ is $k$-stable then $f_{v}$ is $k$-stable. 

\begin{lmma}
Let $v \in E(a)$ and $f \in \mathcal{S}$ be given. Then $\phi(v)^{*}(f+\mathcal{N})=f_{v}+\mathcal{N}$.
\end{lmma}
\textit{Proof.} It suffices to prove that for every $g \in \mathcal{S}$, $\langle \phi(v)^{*}(f+\mathcal{N})|g+\mathcal{N}\rangle=\langle f_{v}+\mathcal{N}|g+\mathcal{N}\rangle$. Let $g \in \mathcal{S}$ be given. Without loss of generality, we can assume that $f$ and $g$ are $k$-stable for some $k \geq 1$. Note that $\phi_{0}(v)g$ is $k+2$-stable. Now calculate as follows to observe that 
\begin{align*}
\langle \phi(v)^{*}(f+\mathcal{N})|g+\mathcal{N} \rangle & = \langle f|\phi_{0}(v)g \rangle \\
                                                                                        &=  \int_{L_{k+2}}\langle f(x)|vg(x-a)\rangle dx ~~(\textrm{by Lemma \ref{inner product1}})\\
                                                                                        &= \int_{L_{k+1}}\langle f(x+a)|vg(x) \rangle dx \\
                                                                                        &= \int_{L_{k+1}}\langle v^{*}f(x+a)|g(x) \rangle dx \\
                                                                                        &=\int_{L_{k+1}}\langle f_{v}(x)|g(x) \rangle dx \\
                                                                                        &=\langle f_{v}+\mathcal{N}|g+\mathcal{N} \rangle ~~(\textrm{by Lemma \ref{inner product1}}). 
 \end{align*}          
This completes the proof. \hfill $\Box$

\begin{thm}
\label{essential}
The representation $\phi$ is essential. 
\end{thm}
\textit{Proof.} Let $\alpha:=\{\alpha_{x}\}_{x \in \Omega}$ be the semigroup of normal $*$-endomorphisms associated to $\phi$. In order to show that $\alpha_{x}$ is unital for each $x \in \Omega$, it suffices to prove that $\alpha_{a}$ is unital. 
Indeed, suppose that $\alpha_{a}$ is unital. Then $\alpha_{na}=\alpha_{a}^{n}$ is unital for every $n \geq 1$. Let $x \in \Omega$ be given. Choose $n \geq 1$ such that $na>x$. Write $na=x+y$ with $y \in \Omega$. Then $1=\alpha_{na}(1)=\alpha_{x}(\alpha_{y}(1)) \leq \alpha_{x}(1)$. Hence $\alpha_{x}$ is unital. Thus it suffices to show that $\alpha_{a}$ is unital. 

Let $\{v_1,v_2,\cdots\}$ be an orthonormal basis for $E(a)$. We claim that \[\sum_{i=1}^{\infty}\phi(v_i)\phi(v_i)^{*}=1\] where the sum is a strongly convergent sum. Since $\{\phi(v_i)\phi(v_i)^{*}\}_{i=1}^{\infty}$ is a sequence of pairwise orthogonal projections, it suffices to show that \[\sum_{i=1}^{\infty}\langle \phi(v_i)\phi(v_i)^{*}(f+\mathcal{N})|f+\mathcal{N}\rangle=||f+\mathcal{N}||^{2}\] for every $f \in \mathcal{S}$. 

Let $f \in \mathcal{S}$ be given. Assume that $f$ is $k$-stable for some $k \geq 1$. Then $f_{v_i}$ is $k$-stable for every $i$. Now calculate as follows to observe that 
\begin{align*}
\sum_{i=1}^{\infty}||\phi(v_i)^{*}(f+\mathcal{N})||^{2}&=\sum_{i=1}^{\infty}||f_{v_i}+\mathcal{N}||^{2} \\
                                                                                 & =\sum_{i=1}^{\infty} \int_{L_{k}}||v_{i}^{*}f(x+a)||^{2}dx ~~(\textrm{by Lemma \ref{inner product1}}) \\
                                                                                 &=\int_{L_{k}}\Big(\sum_{i=1}^{\infty}||v_{i}^{*}f(x+a)||^{2}\Big)dx\\
                                                                                 &=\int_{L_{k}}||f(x+a)||^{2}dx ~~(\textrm{by Equality \ref{onb}})\\
                                                                                 &=\int_{L_{k}}||f(x)||^{2}dx ~~(\textrm{ Since $f$ is $k$-stable})\\
                                                                                 &=||f+\mathcal{N}||^{2}~~(\textrm{by Lemma \ref{inner product1}}).
\end{align*}
In the third equality of the above calculation, we have interchanged the summation and the integral which is permissible since the terms involved are non-negative. This completes the proof. \hfill $\Box$

We end our paper by stating what we have done so far as a Theorem whose proof is a direct consequence of Theorem \ref{essential} and Proposition \ref{Essential}.
\begin{thm}
Let $E:= \coprod_{x \in \Omega}E(x)$ be a product system over $\Omega$ such that $E(x)$ is infinite dimensional for every $x\in \Omega$. Denote the representation of $E$ constructed in Proposition \ref{measurability of phi} by $\phi$. Let $\alpha:=\{\alpha_{x}\}_{x \in P}$ be the $E_0$-semigroup associated to $\phi$. Then $E$ is isomorphic to the product system associated to $\alpha$. 
\end{thm}

\bibliography{references}
 \bibliographystyle{plain}

 \noindent{\sc S.P. Murugan} (\texttt{spmurugan@cmi.ac.in})\\
         {\footnotesize  Chennai Mathematical Institute, \\
Siruseri, 603103, Tamilnadu.}\\[1ex]
{\sc S. Sundar}
(\texttt{sundarsobers@gmail.com})\\
         {\footnotesize  Chennai Mathematical Institute,  \\
Siruseri,  603103, Tamilnadu.}

\end{document}